\documentclass{amsart}
\usepackage{graphicx,amssymb,color}

\newtheorem{theorem}{Theorem}[section]
\newtheorem{lemma}[theorem]{Lemma}
\newtheorem{proposition}[theorem]{Proposition}

\newtheorem{defn}{Definition}[section]
\newtheorem{assumption}{Assumption}

\def \grad{\nabla}

\usepackage{amsmath,amsthm,color,latexsym}   
\usepackage{graphicx}
\usepackage{esint}
\usepackage[hang,flushmargin]{footmisc}
\usepackage{tikz}
\usetikzlibrary{intersections}

\numberwithin{equation}{section}

\title{Inverse Radiative Transport on Non-Convex Domains}

\author[Chung]{Francis J. Chung}
\address{Department of Mathematics, University of Kentucky, Lexington KY 40506, USA}
\email{fj.chung@uky.edu}

\begin{document}

\begin{abstract}
This paper addresses an inverse source problem for the radiative transport equation in non-convex domains with a single ``corner'', generated by the union of two convex domains, asking whether an observer restricted to one part of the boundary of such a domain can ``see around the corner'' and recover a boundary source restricted to a part of the boundary hidden from the observer. The main result shows that with a condition on the geometry of the corner, an analyticity condition on the coefficients of the equation, and a finiteness condition on the scattering kernel, an isotropic boundary source hidden around the corner from the location of the observations can be recovered by measurements of the scattered light.  
\end{abstract}

\maketitle

\section{Introduction}

Light propagation in a scattering medium is an interesting process which gives rise to a number of inverse problems, where optical coefficients or sources must be reconstructed from observations of the scattered light. See for example ~\cite{ArrSch, BalSurvey} for a survey of results where light propagation is modeled by radiative transport, and ~\cite{UhlSurvey, Isa} for a survey of related elliptic problems, which can model light propagation in very highly scattering media ~\cite{ArrSurvey}. In these papers, scattering is the problem that has to be overcome.  But the phenomenon of scattering -- of light changing directions as it travels -- also raises the possibility of seeing things that would be impossible to see otherwise.  

Consider for example the situation of a hallway with a corner, with a light source hidden around the corner. Since for the purpose of this article I'd like to focus on the behaviour of scattered light, let's for the sake of argument ignore the phenomenon of reflection off the walls -- maybe Anish Kapoor painted them all in the ridiculously black colour he's been trying to keep to himself \cite{FinGre}.  

In Kapoor's sombre unreflective hallway, we have no chance of understanding or even perceiving the light source in the absence of scattering -- if light simply travels in straight lines and vanishes when it hits the solid surfaces, we can't see around the corner by definition. But if the hallway is hazy or misty enough that light scatters, then we should be able to perceive light which has traveled around the corner.  The main question of these notes is, can we see around the corner using the scattered light? That is to say, can we reconstruct information about the light source from the observations of the scattered light?

To address this as a mathematical problem, we can use a radiative transport equation (RTE) to model light propagation in a scattering medium.  To do this, let $n \geq 3$, and suppose $\Omega \subset \mathbb{R}^n$ is a piecewise smooth bounded domain. Let $X = \Omega \times S^{n-1} \times [0,\infty)$, and for $u: X \rightarrow \mathbb{R}$, we'll consider the frequency-dependent equation
\begin{equation}\label{FrequencyRTE}
\theta \cdot \grad_x u(x,\theta,v) +\sigma(x,v) u(x,\theta, v)= k(x,v) \int_0^\infty \int_{S^{n-1}}u(x,\theta',v')q(v - v') \, d\theta' dv'.
\end{equation}

Here $u(x,\theta,v)$ is the specific intensity of light at frequency $v$, so $u(x,\theta,v)$ can be thought of as representing the intensity of light at point $x$ traveling in direction $\theta$ with frequency $v$.  Equation \eqref{FrequencyRTE} can be thought of as modeling light propagation in a scattering medium, with the coefficients $k$ and $\sigma$ governing scattering and total attenuation, respectively, and $q$ governing the propensity of incoming light at frequency $v'$ to scatter at frequency $v$.  Note that if we think of $q$ as the delta distribution, then this reduces to the familiar RTE
\begin{equation}\label{SimpleRTE}
\theta \cdot \grad_x u(x,\theta,v) +\sigma(x,v) u(x,\theta, v)= k(x,v) \int_{S^{n-1}}u(x,\theta',v) \, d\theta'
\end{equation}
which acts independently of $v$.   

To discuss the boundary conditions for \eqref{SimpleRTE}, define
\[
\partial \Omega_{\pm} = \{ (x,\theta,v) \in \partial \Omega \times S^{n-1} \times [0,\infty) | \pm \theta \cdot n(x) \geq 0\},
\]
where $n(x)$ represents the outward unit normal at $x$; for a given $x \in \partial \Omega$, we can think of $(x,\theta,v) \in \partial \Omega_+$ as meaning that $\theta$ points out of the domain at $x$, and $(x,\theta,v) \in \partial \Omega_-$ meaning that $\theta$ points into the domain at $x$.  

Now we can distinguish between two kinds of source problems for the RTE: an \emph{interior} source problem, which takes the form
\begin{equation}\label{SVP}
\begin{split}
    \theta \cdot \grad_x u(x,\theta,v) +\sigma(x,v) u(x,\theta, v) &= k(x,v) \int_0^\infty \int_{S^{n-1}}u(x,\theta',v')q(v - v') \, d\theta' dv' \\
    & \, \quad \, +S(x,v)\\
    u|_{\partial \Omega_-} &= 0,
\end{split}
\end{equation}
and a \emph{boundary} source problem, which takes the form 
\begin{equation}\label{BVP}
\begin{split}
    \theta \cdot \grad_x u(x,\theta,v) +\sigma(x,v) u(x,\theta, v) &= k(x,v) \int_0^\infty \int_{S^{n-1}}u(x,\theta',v')q(v - v') \, d\theta' dv'\\
    u|_{\partial \Omega_-}(x,\theta,v) &= f(x,v). \\
\end{split}
\end{equation}

 Note that it follows from results of \cite{BalJol} (see also similar ideas in \cite{BalChuSch}) that the inverse source problem is not solvable if the boundary source $f$ depends on $x$ and $\theta$, hence the restriction to isotropic sources $f$ independent of angle; the interior source $S$ is restricted to be isotropic for the same reason.    

The following proposition ensures that both forward problems are solvable. 

\begin{proposition}\label{Solvability}
Suppose $S, \sigma \in C(\Omega \times [0,\infty))$, $k \in C(\Omega \times [0,\infty))$, $f \in C(\partial \Omega \times [0,\infty))$, $q \in C'(\mathbb{R})$ normalized so that $\|q\|_{C'} = 1$, and 
\[
\sigma(x,v) \geq k(x,v) \geq 0 
\]
for all $x \in \Omega$.  The equation 
\begin{equation}\label{DSVP}
\begin{split}
    \theta \cdot \grad_x u(x,\theta,v) +\sigma(x,v) u(x,\theta,v) &= k(x,v) \int_0^\infty \int_{S^{n-1}}u(x,\theta',v')q(v - v') \, d\theta' dv' \\
    &\, \quad \, +S(x,v)\\
    u|_{\partial \Omega_-} &= f(x,v)
\end{split}
\end{equation}
has a unique solution $u \in C(X)$.
\end{proposition}

In all subsequent results in these notes, we will assume the hypotheses of Proposition \ref{Solvability}. Now we turn to the inverse problem. Suppose $E \subset \partial \Omega$ is an open set.  We define 
\[
E_+ = \{(x,\theta,v) \in \partial \Omega_+| x \in E\}.
\]
Each of \eqref{SVP} and \eqref{BVP} now yields an inverse problem: does $u|_{E_+}$ determine the source function?

If $E$ covers the entire boundary, then this problem has been extensively investigated in the frequency independent case $q(v) = \delta(v)$ (see \cite{BalTam, Lar, Sie, SteUhl}, plus hybrid versions of the problem in \cite{BalChuSch, ChuYanYan}). In the case where $E$ is a strict subset of $\partial \Omega$, again in the frequency independent case, the inverse interior source problem defined by \eqref{SVP} was investigated by Hubenthal ~\cite{Hub} (see also ~\cite{FujSadTam} for a similar result in two dimensions) who concluded that $S$ can be recovered if its support a priori lies within the set of $x$ for which every line through $x$ intersects $E$.  Most of these results also cover more general scattering kernels $k$, allowing anisotropic scattering.  

In order to do better, we will make some analyticity assumptions on $\sigma$ and $k$.  
\vspace{2mm}

\begin{assumption}\label{AnalAssumption}[Analyticity Assumption]
Assume that $\sigma(x,v)$ is constant as a function of $x$ for each fixed $v$, and $k(x,v) > 0$ is analytic as a function of $x$ on a neighbourhood of $\Omega$, uniformly in $v$.
\end{assumption}

With this assumption in place, we can obtain the following result for convex domains.  

\begin{proposition}\label{ConvexProp}
Suppose $\Omega$ is convex, and $E$ is an open subset of $\partial \Omega$.  Under the conditions of assumption \ref{AnalAssumption}, if $u$ solves \eqref{SVP} and $u|_{E_+} = 0$, then $S \equiv 0$.
\end{proposition}

But the real point of these notes is to consider non-convex domains. For these domains we will confine ourselves to considering the simpler inverse boundary source problems defined by \eqref{BVP}.  We will also need an additional set of assumptions on $q$ and $f$. 
\vspace{2mm}

\begin{assumption}\label{RamanAssumption}[Raman Scattering Assumption]
Assume that $q$ is supported only in $(-\infty,-c]$ for some fixed positive constant $c$, and that $f(x,v)$ is supported only for $0< v \leq M$, for some fixed upper bound $M > c$. We also assume that the operator $g \mapsto \int_0^{\infty}g(v') q(v-v')\, dv'$ is partially invertible in the sense that $\{\int_0^{\infty}g(v') q(v-v') dv'| v > 0\}$ can be used to recover $g(v_0)$ for some fixed $0<v_0 \leq M$.

\end{assumption}

This is meant to simulate a Raman scattering process \cite{HubMih}, where frequency (equivalently, energy) of scattered photons decreases after scattering \cite{HubMih}.  The last assumption can be shown for a range of $q$, e.g.  $q(v) = \delta(v + c)$. Note that the combination of Raman scattering and isotropic scattering (where $k$ is independent of angle) is somewhat unusual, but see also further discussion in section 4.  

Finally, we'll need one more set of assumptions to constrain the type of non-convexity in our domain, which brings us to the following definition.

\begin{defn}\label{Recoverable}
    Suppose $\Omega$ is a piecewise smooth bounded domain, $E \subset \partial \Omega$ is open, and $V \subset \Omega$.  We define $V$ to be \emph{recoverable} from $E$ if every point in $V$ lies in a plane $P$ such that $P\cap E$ is nonempty and $P \cap \Omega$ is convex.  
\end{defn}

\vspace{4mm}

\begin{center}

\includegraphics[scale=0.1]{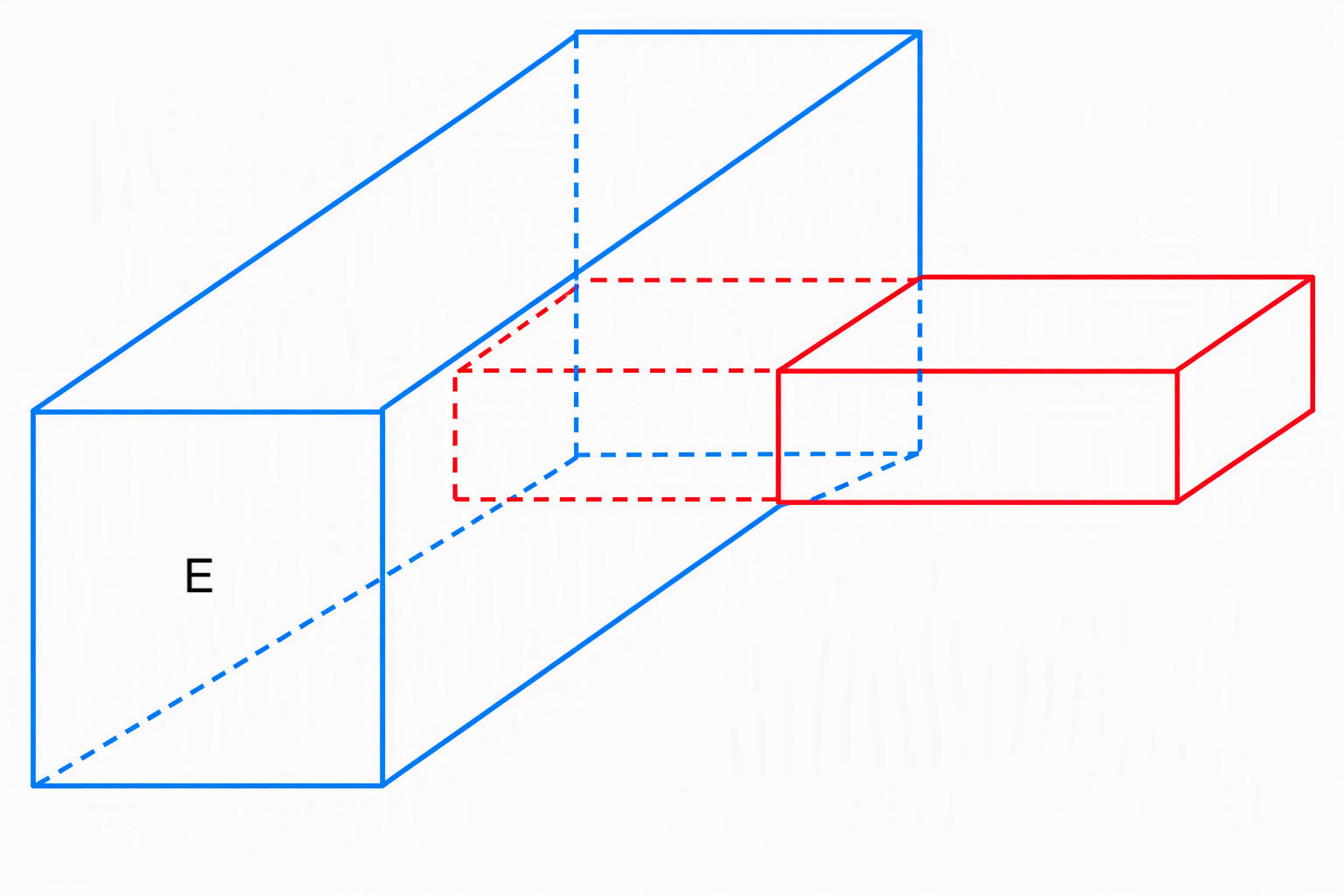}

{\small Fig. 1. An example of a recoverable corner: $\Omega$ is the union of the red and blue rectangular prisms.  The entire (solid) blue prism is recoverable from the labeled front face $E$; in particular the entire intersection of the two prisms is recoverable from $E$.  }
\end{center}

\vspace{4mm}

\begin{center}

\includegraphics[scale=0.1]{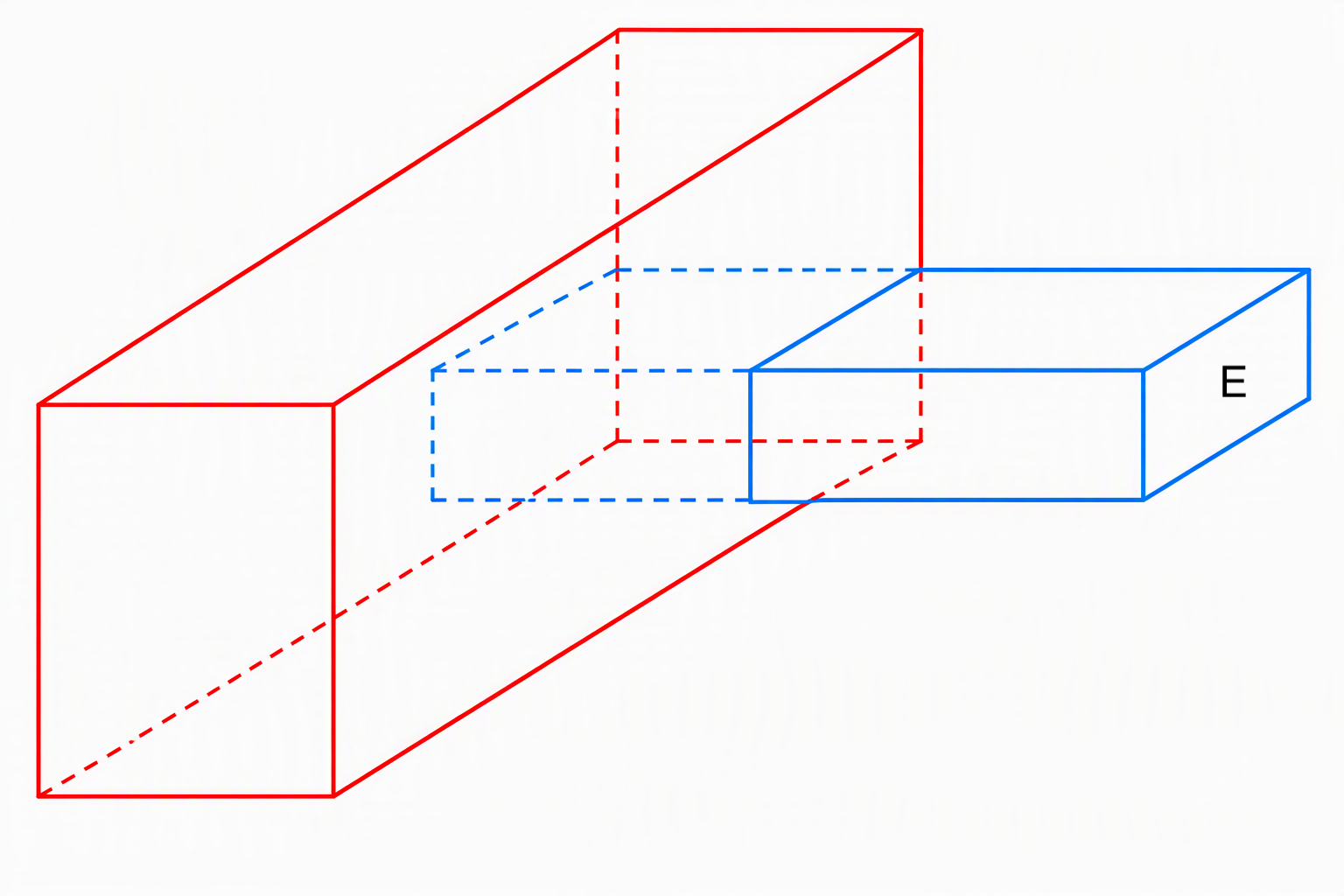}

{\small Fig. 2. An example of an unrecoverable corner.  Again $\Omega$ is the union of the red and blue rectangular prisms, but in this case no open subset of the overlap is recoverable from $E$, the labeled rightmost face of the blue prism.  }
\end{center}

\vspace{4mm}

With all of this in place, we can state the main result of this article. 

\begin{theorem}\label{MainTheorem}
Suppose $\Omega = A \cup B$ with $A \cap B$ having nonempty interior, and $A$ and $B$ piecewise smooth, bounded, and convex. Let $E \subset \partial \Omega$ be an open set such that $B$ is recoverable from $E$.  If  $f \in C(\partial \Omega \times [0,M])$ is supported only for $x \in F \subset \partial A$ which is not connected to $E$ by straight lines through $\Omega$, and $u$ solves \eqref{BVP} under assumptions \ref{AnalAssumption} and \ref{RamanAssumption}, then $u|_{E_+} = 0$ implies that $f(x,v_0)= 0$ for all $x \in F$.
\end{theorem}

\vspace{4mm}
\begin{center}

\includegraphics[scale=0.1]{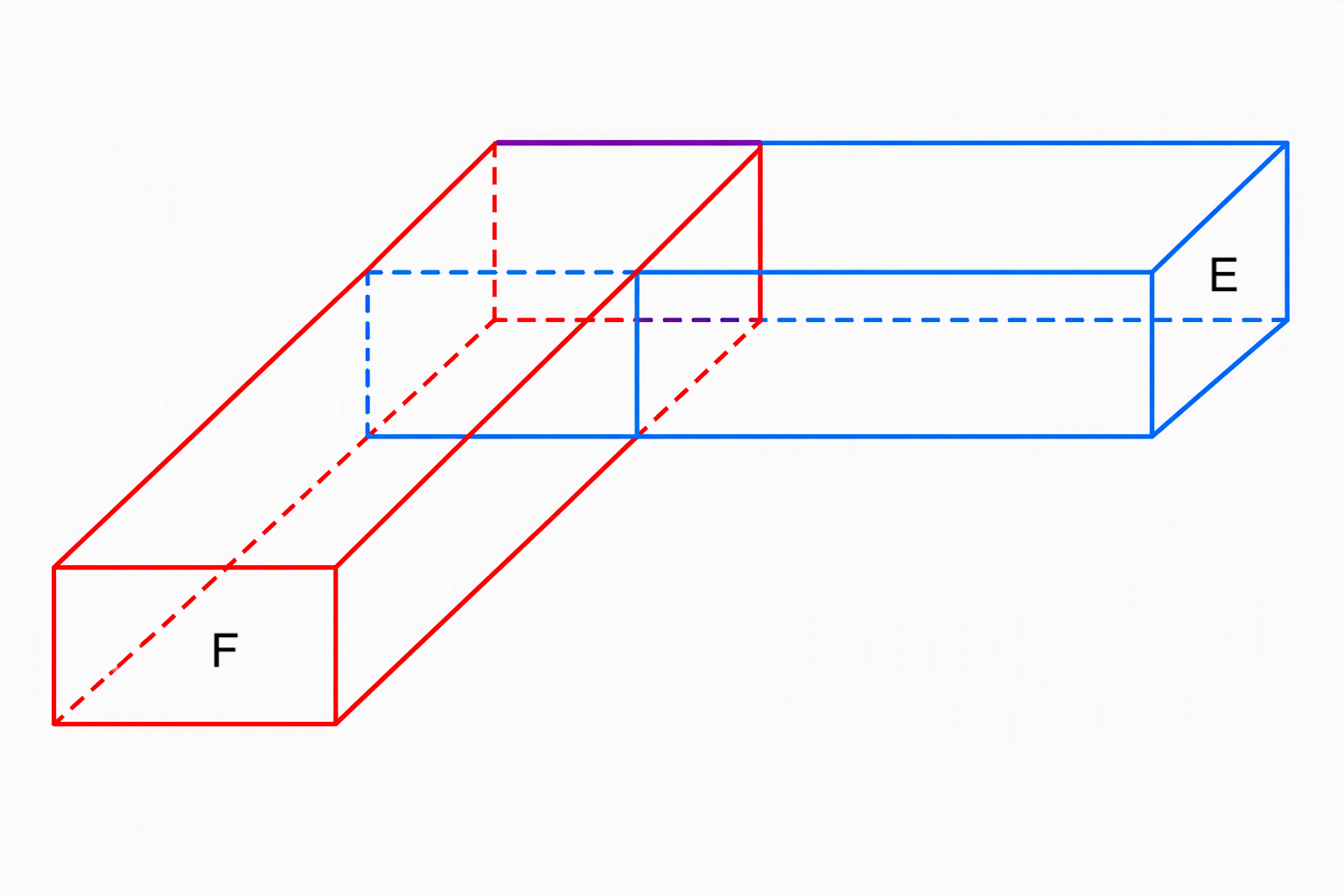}

{\small Fig. 3. An example of the situation covered by Theorem \ref{MainTheorem}.  $A$ and $B$ are the red and blue sets, respectively, with $F$ and $E$ corresponding to the labeled faces.  Theorem \ref{MainTheorem} tells us that isotropic sources restricted to $F$ can be recovered by measurements on $E$. }
\end{center}
\vspace{4mm}

The remainder of these notes are organized as follows.  In Section 2 we will introduce some notation and preliminary results and prove Proposition \ref{ConvexProp}.  Section 3 is devoted to the proof of Theorem \ref{MainTheorem}.  Finally, Section 4 contains discussion related to the necessity of the various conditions in the statement of Theorem \ref{MainTheorem}.
\vspace{4mm}

\noindent \textbf{Acknowledgements.} The author wishes to thank Todd Quinto, Gaik Ambartsoumian, and Venky Krishnan for the conference invitation that prompted the finishing of this work, as well as Teemu Saksala for his urging to actually finish a draft of the manuscript and make it public.  This work was partially supported by Simons collaboration grant 582020.
\vspace{4mm}

Figures 1-3 were made with the assistance of Microsoft Copilot. In addition, several AI products were consulted during this project to look for references for various results. In the case of Lemma 3.1, no reference was located, but ChatGPT suggested a method of proof which is the basis of the one given in this paper.  A similar situation occurred with the analysis of layer potentials on page 14 (although here a reference was suggested and consulted). The remainder of this paper was generated by the author in the traditional fashion. All errors and omissions are the sole responsibility of the author.

\section{Preliminaries}

First we'll need some notation to describe solutions to \eqref{BVP}. Suppose $x, y \in \bar{\Omega}$, such that the line segment from $x$ to $y$ is contained in $\Omega$. We define the attenuation factor
\begin{equation}\label{Attenuation}
\alpha_v(x,y) = \exp\left(- \int_0^{\|y-x\|}\sigma(x - s\theta,v) \, ds\right)
\end{equation}
where $\theta$ is the unit vector in the direction of the vector $x-y$. Note that 
\[
\alpha_v(x,y) = \alpha_v(y,x).
\]

Using this notation we find that the solution to the transport equation
\begin{equation}\label{TEquation}
\begin{split}
\theta \cdot \grad u + \sigma(x,v) u &= Q(x,v)\\
u|_{\partial \Omega_-} &= 0 \\
\end{split}
\end{equation}
is given by 
\begin{equation}\label{TInverseDefn}
u(x,\theta,v) = (T^{-1}Q)(x,\theta,v) := \int_0^{\|x - x_{\theta-}\|} \alpha_v(x, x - t\theta)Q(x - t\theta,v)dt.
\end{equation}
where $x_{\theta-}$ is the first point on $\partial \Omega$ obtained by traveling in the $-\theta$ direction from $x$.  

Note that since $\alpha_v \leq 1$ and $\|x - x_{\theta-}\|$ is uniformly bounded on $\overline{\Omega}$, 
\begin{equation}\label{TInverseBound}
\|T^{-1} Q\|_{C(X)} \lesssim \|Q\|_{C(\Omega \times [0,\infty))}.
\end{equation}

 Using this notation we can write a solution to \eqref{DSVP} explicitly in terms of an infinite series.  See also similar proofs in \cite{Chu,EggSch}.  

\begin{proof}[Proof of Proposition \ref{Solvability}]

First let's address the problem \eqref{DSVP} in the case that $f \equiv 0$; in other words, when $u$ solves \eqref{SVP}.  Since the RTE in \eqref{SimpleRTE} takes the form of \eqref{TEquation} with 
\[
Q(x,v) = k(x,v) \int_0^\infty \int_{S^{n-1}}u(x,\theta',v')q(v - v') \, d\theta' dv' +S(x,v),
\]
it follows from \eqref{TInverseDefn} that $u$ solves \eqref{SVP} when 
\begin{equation}\label{BigUExpression}
u(x,\theta,v) = T^{-1}\left(k(x,v) \int_0^\infty \int_{S^{n-1}}u(x,\theta',v')q(v - v') \, d\theta' dv' +S(x,v)\right).
\end{equation}

Let 
\begin{equation}\label{KDefn}
Ku(x,v) = k(x,v) \int_0^\infty \int_{S^{n-1}}u(x,\theta',v')q(v - v') \, d\theta' dv',
\end{equation}
so we can rewrite \eqref{BigUExpression} as 
\begin{equation}\label{SimpleUExpression}
u = T^{-1}(Ku + S)
\end{equation}
or
\[
u -T^{-1}Ku = T^{-1}S.
\]
Direct calculation combined with the inequality  $\sigma \geq k$ and the normalization of $q$ gives for some $c < 1$ depending on $\Omega$ the result that 
\[
\|T^{-1}Ku\|_{C(X)} \leq c\|u\|_{C(X)} 
\]
for all $u \in C(X)$.  It follows that the infinite series 
\begin{equation}\label{CollisionExpansion}
u = \sum_{j=0}^{\infty}(T^{-1}K)^{j}T^{-1}S
\end{equation}
converges and gives the unique solution to \eqref{SimpleUExpression}.

To expand this proof to the case when $f \neq 0$, we only need a short additional argument.

For $f \in C(\partial \Omega \times [0,M])$, define $Jf$ on $X$ by 
\begin{equation}\label{JDefn}
Jf(x,\theta,v) = \alpha_v(x,x_{\theta-})f(x_{\theta_-},v).
\end{equation}

Now by the proof for the solvability of \eqref{SVP}, the problem 
\begin{equation*}
\begin{split}
    \theta \cdot \grad_x w +\sigma w &= Kw + KJf\\
    w|_{\partial \Omega_-} &= 0
\end{split}
\end{equation*}
has a unique solution $w$, and one can check that 
\begin{equation}\label{BoundaryU}
u = w + Jf 
\end{equation}
solves \eqref{BVP}. To obtain a solution to \eqref{DSVP}, we can add the solutions to \eqref{SVP} and \eqref{BVP}.  

\end{proof}

The proofs of Proposition \ref{ConvexProp} and Theorem \ref{MainTheorem} rely on a limited data result for the attenuated X-ray transform.  For this we'll adapt a proposition from \cite{Qui93c} (see Theorem 1.1 and the discussion preceding its statement;  also see \cite{BomQui, BomQui93, Qui94}).

To set up the proposition, let $Y$ denote the space of lines in $\mathbb{R}^2$, and let $\mu: \mathbb{R}^2 \times Y \rightarrow \mathbb{R}$.  For a function $f \in C(\Omega)$, define for each $y \in Y$
\[
R_{\mu}f(y) = \int_{y} f(x) \mu(x,y) \, dm(x)
\]
where $dm$ is the usual measure on the line $y \subset \mathbb{R}^2$.

\begin{proposition}\label{AttenuatedPartialXRay}
Let $f \in C(V)$ for some compact set $V \subset \mathbb{R}^2$, and suppose $\mu: \mathbb{R}^2 \times Y \rightarrow \mathbb{R}$ is analytic and non-zero on $V$.  Let $\mathcal{A}$ be an open connected subset of $Y$ such that there exists at least one $y_0 \in \mathcal{A}$ with $f =0$ on a neighbourhood of $y_0$.  If $R_{\mu}f(y) = 0$ for all $y \in \mathcal{A}$, then $f \equiv 0$ on each $y \in \mathcal{A}$.  \cite{Qui93c}
\end{proposition} 

With this in hand, we can prove Proposition \ref{ConvexProp}.

\begin{proof}[Proof of Proposition \ref{ConvexProp}]

Let $P$ be a plane that intersects $E$, and consider for the moment the data confined to $P\cap E$. More precisely, let $(x,\theta,v) \in (P\cap E)_+$ with $\theta$ parallel to $P$. Equation \eqref{SimpleUExpression} implies that 
\[
u(x,\theta,v) = T^{-1}(Ku + S)(x,\theta,v).
\]
Fixing a $v \in [0,\infty)$ and rewriting the definition of $T^{-1}$ in \eqref{TInverseDefn}, we get  
\[
T^{-1}Q(x,\theta,v) = \int_y \alpha_v(x,z) Q(z,v) dm(z) 
\]
where $y$ is the line in $P$ that goes through $x$ in direction $\theta$.  By the hypotheses of the proposition, it follows that
\[
\int_y \alpha_v(x,z) (Ku + S)(z) dm(z)  = 0.
\]
for each line $y \subset P$ that intersects $E$.  We want to use this to show that $R_{\mu} (Ku+S)(y) = 0$ for some analytic weight $\mu$.  To do this, choose coordinates on $P$, and for each $\omega \in [0,\pi]$, let
\begin{align*}
\vec{\omega} &= (\cos \omega, \sin \omega) \in P\\
\vec{\omega}^{\perp} &= (\cos(\omega + \pi/2), \sin(\omega + \pi/2)) \in P. \\ 
\end{align*}
Now we can parametrize $Y$ as $[0,\pi] \times \mathbb{R}$ by writing each $y \in Y$ as $y = y(\omega, p) = \{x \in P| x \cdot \vec{\omega} = p\}$.   Since $\Omega$ is convex, Proposition 3.2 in \cite{Chu} implies that there exists an open set $V \subset P$ disjoint from $P \cap \Omega$ such that for every line $y$ in $P$ that intersects $V$, $y$ intersects $E$. Now we can pick coordinates to ensure that all lines through $V$ have $x \in y \cap \Omega$ only if $x \cdot \vec{\omega}^{\perp} > 0$.

Then set 
\[
\mu(x,y) = \mu(x,\omega, p) = e^{-\int_0^1\sigma_v(x\cdot \vec{\omega}^{\perp}\vec{\omega}^{\perp}t + p\vec{\omega},v)x\cdot \vec{\omega}^{\perp} \, dt}
\]
and note that for $x \in y \cap \Omega$, we have 
\begin{align*}
\mu(x,y) &= e^{-\int_0^1\sigma_v(|x\cdot \vec{\omega}^{\perp}|\vec{\omega}t + p\omega,v)|x\cdot \vec{\omega}^{\perp}| \, dt} \\
&= e^{-\int_0^1\sigma_v(|x-p\vec{\omega}|\vec{\omega}t + p\omega)|x- p\vec{\omega}| \, dt}. \\
\end{align*}
By a change of variables, this gives us
\[
\mu(x,y) = \alpha_v(x, p\vec{\omega}).
\]

Then since $p\vec{\omega} \cdot \vec{\omega}^{\perp} = 0$, it follows that $p\vec{\omega} \notin \Omega$, so for $z \in \Omega$ and $y$ a line through $z$ in $P$ intersecting $E$ at $x$, we must have $x$ lying between $p\vec{\omega}$ and $z$ in $y$.  Therefore
\begin{align*}
R{\mu}(Ku+S)(y) &= \int_y \mu(z,y) (Ku+S)(z)dm(z) \\ 
			    &= \int_y \alpha_v(z,p\vec{\omega}) (Ku+S)(z)dm(z) \\ 
			    &= \int_y \alpha_v(x,z)\alpha_v(x,p\vec{\omega}) (Ku+S)(z)dm(z) \\ 
			    &= \alpha_v(x,p\vec{\omega})\int_y \alpha_v(x,z) (Ku+S)(z)dm(z) \\ 			    
			    &= 0 \\ 			    			    
\end{align*}
for each line $y$ in $P$ that intersects $E$, and hence for any line $y$ in $P$ that intersects $V$.

Note that since $\sigma(\cdot,v)$ is constant for each fixed $v$, $\mu$ is analytic and non-zero. The set $\mathcal{A}$ of lines in $P$ that intersect $V$ is open, and moreover since $P\cap \Omega$ is convex, $\mathcal{A}$ contains at least one line with positive distance from the support of $Ku + S$.  It now follows from Proposition \ref{AttenuatedPartialXRay}  that $Ku + S = 0$ on $P \cap \Omega$.  Repeating with each plane $P$ that intersects $E$, we find that $Ku + S \equiv 0$ on $\Omega$.  Then equation \eqref{SimpleRTE} says that 
\[
\theta \cdot \grad u + \sigma u = 0
\]
on $X$, with boundary condition $u|_{\partial \Omega_-} = 0$.  Therefore $u \equiv 0$ on $X$ and therefore $S \equiv 0$ as well. 
\end{proof}

\section{Seeing Around Corners}

First we'll need a pair of lemmas regarding analyticity. The first is a general statement regarding the analyticity of integral transforms. 

\begin{lemma}\label{AnalyticIntegral}
Suppose $K \subset \mathbb{R}^n$ is compact, and $Q \in C(K)'$, the dual space to continuous functions on $K$. Let $P: \mathbb{R}^n \times \mathbb{R}^n \rightarrow \mathbb{R}$ be continuous, and suppose that for some fixed $x_0 \in \mathbb{R}^n$, $P(x,y)$ has a power series expansion
\[
P(x,y) = \sum_{\alpha} a_\alpha(y)(x-x_0)^{\alpha}
\]
which converges absolutely and uniformly for all $(x,y) \in \overline{B_r(x_0)} \times K$.  Then
\[
\psi(x) = \int_{\mathbb{R}^n} P(x,y)Q(y) \, dy
\]   
is analytic and equal to its power series representation in $B_r(x_0)$.
\end{lemma}

\begin{proof}
By the hypotheses on $P$, it follows that for each $y \in K$ and $x \in B_r(x_0)$, we can write
\[
\psi(x) = \int_{\mathbb{R}^n} \sum_{\alpha} a_\alpha(y)(x-x_0)^{\alpha} Q(y) \, dy.
\]
Since the convergence is uniform and $Q$ is supported on a compact set, we can interchange the sum and integral and write
\begin{equation}\label{PsiPower}
\psi(x) = \sum_{\alpha} \left( \int_{\mathbb{R}^n}  a_\alpha(y) Q(y) \, dy \right)(x-x_0)^{\alpha}.
\end{equation}
This gives an expression for $\psi$ as a power series centred at $x_0$.  To check the radius of convergence of this series, note that for $\|x-x_0\| < r$, we have
\[
\sum_{\alpha} \left| \int_{\mathbb{R}^n}  a_\alpha(y) Q(y) \, dy \right| \|x-x_0\|^{\alpha} 
\leq \sum_{\alpha} \int_{\mathbb{R}^n}  |a_\alpha(y)|\|x-x_0\|^{\alpha} |Q(y)| \, dy.  
\]
Since the absolute convergence of the series is also uniform, we can interchange the sum and integral once again and obtain
\[
 \int_{\mathbb{R}^n}  \left( \sum_{\alpha} |a_\alpha(y)|\|x-x_0\|^{\alpha}  \right)|Q(y)| \, dy \leq \int_{\mathbb{R}^n} C |Q(y)| \, dy = \int_{K} C |Q(y)| \, dy
\]
where $C$ is the maximum given by the absolute convergence of the power series for $P$ on the compact set $\overline{B_r(x_0)} \times K$.  Since the integral on the right is finite, it follows that the power series \eqref{PsiPower} converges with radius of convergence $r$.  
\end{proof}

The second lemma is a statement about the analyticity of solutions to the RTE \eqref{BVP} with analytic coefficients.

\begin{lemma}\label{AnalyticRTE}
Suppose $U$ is a convex piecewise smooth bounded domain and $f \in C(\partial U \times [0,\infty))$. If $u$ is the solution to \eqref{BVP}, then for each fixed $v$, $\langle u \rangle(x,v)$  is analytic as a function of $x$ in $U$. 
\end{lemma}

\begin{proof}
It follows from Proposition \ref{Solvability} that 
\begin{equation}\label{UExpansion}
\langle u \rangle = \langle Jf \rangle + \sum_{j=1}^{\infty} \langle (T^{-1}K)^{j} Jf \rangle.
\end{equation}

Note that given Assumption \ref{RamanAssumption}, $f(x,v)$ and hence $Jf(x,v)$ are only supported for $v < M$ for some constant $M$. Moreover, since $q$ is only supported for $v < -c < 0$, for some constant $c$, we get that $KJf$ and hence $T^{-1}KJf$ are supported only for $v < M-c$. After finitely many steps, it follows that $(T^{-1}K)^{j}Jf = 0$ for all $v \geq 0$, and hence the infinite sum above is actually finite:
\begin{equation}\label{FiniteUExpansion}
\langle u \rangle = \langle Jf \rangle + \sum_{j=1}^{N} \langle (T^{-1}K)^{j} Jf \rangle.
\end{equation}
for some $N$.  

To show that $\langle u \rangle$ is analytic, it thus suffices to check each term in the series.  Let's begin by analyzing $Jf$.   Let $x \in U$; from \eqref{JDefn} we have
\[
\langle Jf\rangle(x,v)  = \int_{S^{n-1}} \alpha_v(x,x_{\theta-})f(x_{\theta_-},v) \, d\theta.
\]
Since $U$ is convex, each $y \in \partial U$ takes the form of $x_{\theta-}$ for exactly one $\theta$.  Therefore we can rewrite this as an integral over $\partial U$:
\begin{equation}\label{JConvo1}
\langle Jf\rangle(x,v)  = \int_{\partial U} \alpha_v(x,y)f(y,v)\frac{\vec{n}(y)\cdot(y-x)}{|y-x|^n} \, dS(y),
\end{equation}
where $\vec{n}$ is the outward normal on $\partial U$.  Expanding the dot product, we can write this as a sum of $n$ integrals of the form
\[
\int_{\partial U} A_i(x,y)\phi_i(y,v) \, dS(y).
\]
Each of these integrals can be written as the integral over $\mathbb{R}^n$ by replacing $\phi_j$ with a distribution $\Phi_j(y,v)$ supported on $\partial U$. This gives
\begin{equation}\label{JConvo}
\langle Jf\rangle(x,v) = \sum_{i=1}^n \int_{\mathbb{R}^n}  A_i(x,y)\Phi_i(y,v) \, dy
\end{equation}

Since $\sigma$ is constant, it follows that 
\begin{equation}\label{ADefn}
A_i(x,y)  = \frac{\alpha_v(x,y)(x-y)_i}{\|x-y\|^{n}}
\end{equation}
is analytic in $x$ and $y$ for $x \neq y$, uniformly in $v$. Since $\Phi_i(y,v)$ is only supported for $y \in \partial U$, Lemma \ref{AnalyticIntegral} applies, and tell us that for all $x \in U$, each of these terms is analytic at $x$ with the same radius of convergence as the $A_i$. Therefore in particular this analyticity is uniform in $v$, and the sum $\langle Jf \rangle$, is analytic at $x$ as well, with the same uniformity in $v$.

The next claim is that $\langle (T^{-1}K)^{j} Jf \rangle$ is analytic as well.  Note that 
\[
\langle (T^{-1}K)^j Jf \rangle(x,v) = \langle T^{-1}k(x) \int_0^{\infty} \langle (T^{-1}K)^{j-1}Jf \rangle(x,v') q(v-v') \,dv'
\]
Note that if $\langle (T^{-1}K)^{j-1}Jf \rangle(x,v)$ is analytic as a function of $x$, uniformly in $v$, then an argument like the one in Lemma 3.1 shows that 
\[
\int_0^{\infty} \langle (T^{-1}K)^{j-1}Jf \rangle(x,v') q(v-v') \,dv'
\]
is analytic as a function of $x$ as well, uniformly in $v$.  Since $k(x,v)$ is analytic in $x$, uniformly in $v$, it therefore suffices to show that $\langle T^{-1}Q \rangle(x,v)$ is analytic for $x$ in $U$, uniformly in $v$, if $Q(x,v)$ is. So suppose for fixed $v$, $Q$ is analytic in $U$ as a function of $x$, and write
\[
\langle T^{-1} Q \rangle(x,v) = \int_{S^{n-1}} \int_0^{\|x - x_{\theta -}\|} \alpha_v(x,x-t\theta) Q(x-t\theta,v) \, dt \, d\theta.
\]
Since $U$ is convex, this is equivalent to 
\[
\langle T^{-1} Q \rangle(x,v) = \int_{S^{n-1}} \int_0^{\infty} \alpha_v(x,x-t\theta) Q(x-t\theta,v)\chi_U(x - t\theta) \, dt \, d\theta.
\]
where $\chi_U$ is the characteristic function of $U$.  By understanding this as an integral in spherical coordinates centred at $x$, we obtain
\[
\langle T^{-1} Q \rangle(x,v) = \int_{\mathbb{R}^n} \frac{\alpha_v(x, y)}{\|x-y\|^{n-1}} Q(y,v)\chi_U(y) \, dy.
\]
If $\sigma_v$ is constant then 
\[
\alpha_v(x,y) = \exp\left(- \int_0^{\|y-x\|}\sigma_v \, ds\right) = e^{-\sigma_v\|x-y\|},
\]
so if we write 
\begin{equation}\label{AvDefn}
A_v(x) = \frac{e^{-\sigma_v \|x\|}}{\|x\|^{n-1}}
\end{equation}
then we have 
\[
\langle T^{-1}Q \rangle(x,v) = \int_{\mathbb{R}^n} A_v(x-y)Q(y,v)\chi_U(y) \, dy.
\]
For each fixed $v$ this can be understood as the convolution of $A_v$ with $Q \chi_U$, neither of which are analytic everywhere. However, $A_v$ is analytic away from $0$, so we can break this integral into two parts, using a smooth function $\phi_B$ which is identically one on a fixed ball $B$ containing $x$, with support compactly contained in $U$.  Then 
\[
\langle T^{-1}Q \rangle(x,v) = \int_{\mathbb{R}^n} A_v(x-y)Q(y,v)\chi_{U} (y)(1-\phi_B)(y) \, dy + \int_{\mathbb{R}^n} A_v(x-y)Q(y,v)\phi_{B} (y) \, dy.
\]

Now for $x \in B$, the first integral is the convolution of a function $Q(\cdot,v)\chi_{U}(1-\phi_B)$ with a function $A_v$ which is analytic on the support of $Q(\cdot,v)\chi_{U}(1-\phi_B)$, and hence this is analytic by Lemma \ref{AnalyticIntegral}.  For the second integral, we can apply derivative estimates directly.

To do this, write
\[
Z(x,v) =  \int_{\mathbb{R}^n} A_v(x-y)Q(y,v)\phi_{B} (y) \, dy =  \int_{\mathbb{R}^n} A_v(y)Q(x-y,v)\phi_{B} (x-y) \, dy,
\]
so
\[
\partial^{\alpha}_xZ(x,v) = \int_{\mathbb{R}^n} A_v(y)\partial^{\alpha}_x(Q(x-y,v)\phi_{B} (x-y)) \, dy.
\]
Using the product rule, we get
\begin{equation}\label{FullProductRule}
\partial^{\alpha}_xZ(x,v) = \sum_{\beta \leq \alpha} \binom{\alpha}{\beta} \int_{\mathbb{R}^n}A_v(y)\partial_x^{\alpha - \beta} Q(x-y,v) \partial_x^{\beta} \phi_B(x-y) \, dy.
\end{equation}
We can bound each term in the sum individually. If $|\beta| \leq 1$, we can use the analytic derivative estimates on $Q$, which is analytic in the support of $\phi_B$, to conclude that 
\begin{align*}
& \binom{\alpha}{\beta} \left| \int_{\mathbb{R}^n}A_v(y)\partial_x^{\alpha - \beta} Q(x-y,v) \partial_x^{\beta} \phi_B(x-y) \, dy \right| \\
\leq \, &\binom{\alpha}{\beta} C_{Q,B}(\alpha - \beta)!R_{Q,B}^{-|\alpha - \beta|}\int_{\mathbb{R}^n}| A_v(y) \partial_x^{\beta}\phi_B(x-y)| dy \\
\leq \, &\binom{\alpha}{\beta} C_{Q,B}(\alpha - \beta)!R_{Q,B}^{-|\alpha - \beta|}C_{A,B} \\
\end{align*}
where $C_{Q,B}$ and $R_{Q,B}$ are constants depending on $Q$ and $B$, and $C_{A,B}$ is a constant depending on $A$ and $B$. (Note that $C_{A,B}$ depends directly on the integral of the first derivatives of $\phi_B$, but these depend only on $\phi_B$, which is itself only constrained by the choice of $B$.)

Unfortunately this strategy does not work for all the terms, because $\phi_B$ is not analytic and hence we have no bound on arbitrarily high derivatives of $\phi_B$.  But for $|\beta| > 1$, we can transfer all but one derivative from $\phi$ back onto the remaining functions in the integrand: if $|\beta| > 1$ then 
\begin{align*}
& \binom{\alpha}{\beta} \left| \int_{\mathbb{R}^n}A_v(y)\partial_x^{\alpha - \beta} Q(x-y,v) \partial_x^{\beta} \phi_B(x-y) \, dy \right| \\
= \, &\binom{\alpha}{\beta} \left| \int_{\mathbb{R}^n} \partial_y^{\beta'}(A_v(y) \partial_y^{\alpha - \beta} Q(x-y,v))\partial_{x_i}\phi_B(x-y) \, dy \right|
\end{align*}
where we've rewritten $\partial^{\beta}_{x} = \partial^{\beta'}_x\partial_{x_i}$ for some first order derivative operator $\partial_{x_i}$; note that the derivatives in $x$ and $y$ are equivalent up to sign, because of the arguments of the form $x-y$ in the relevant functions. Another application of the product rule gives us 
\[
\binom{\alpha}{\beta} \left| \sum_{\gamma \leq \beta'} \binom{\beta'}{\gamma}\int_{\mathbb{R}^n}  \partial_y^{\gamma} A_v(y) \partial_y^{\alpha - \beta + \beta' - \gamma} Q(x-y,v) \partial_{x_i}\phi_B(x-y)\, dy \right|.
\]

Since $\phi_B \equiv 1$ on $B$, it follows that $A(y)$ is also analytic on the support of $\partial_{x_i} \phi_B(x-y)$, so we can take advantage of the analytic derivative estimates on $A$ and $Q$.  Doing this gives us the bound 
\[
\sum_{\gamma \leq \beta'}\binom{\alpha}{\beta}\binom{\beta'}{\gamma}C_{A,B} \gamma ! R_{A,B}^{-|\gamma|}C_{Q,B}(\alpha - \beta +\beta' - \gamma)! R_{Q,B}^{-|\alpha - \beta+\beta' -\gamma|} \int_{\mathbb{R}^n} |\partial_{x_i}\phi_B(x-y)| \, dy,
\]
and since the integral then only depends on the choice of $\phi_B$, which is only restricted by the choice of $B$, we can combine constants and write 
\[
\sum_{\gamma \leq \beta'}\binom{\alpha}{\beta}\binom{\beta'}{\gamma}C_{A,Q,B} \gamma! (\alpha - \beta +\beta' - \gamma)! R_{A,B,Q}^{-|\alpha -1|}.
\]
Expanding the chooses, we get 
\[
\sum_{\gamma \leq \beta'}\frac{\alpha !}{(\alpha - \beta)! \beta!}\frac{\beta'!}{\gamma! (\beta' - \gamma)!}C_{A,Q,B} \gamma! (\alpha - \beta +\beta' - \gamma)! R_{A,B,Q}^{-|\alpha -1|},
\]
which is bounded above by 
\[
\sum_{\gamma \leq \beta'}C_{A,Q,B}   R_{A,B,Q}^{-|\alpha -1|}\alpha! \binom{\alpha - \beta +\beta' - \gamma}{\alpha - \beta},
\]
which in turn is bounded above by 
\[
C_{A,Q,B}   R_{A,B,Q}^{-|\alpha -1|}\alpha! \binom{\alpha}{\alpha - \beta} |\alpha|.
\]
Combining this estimate with the estimate on the terms with $|\beta| \leq 1$, we get that each term in \eqref{FullProductRule} is bounded above by a term of the form 
\[
C_{A,Q,B}   R_{A,B,Q}^{-|\alpha|}\alpha! \binom{\alpha}{\beta}
\]
for some choice of constants $C_{A,Q,B}$ and $R_{A,B,Q}$ depending only on $A,Q,$ and $B$.  Summing over all terms then gives an upper bound of 
\[
|\partial^{\alpha}_xR(x,v) | \leq C_{A,Q,B}   \left(\frac{R_{A,B,Q}}{3}\right)^{-|\alpha|}\alpha!
\]
which lets us conclude that $R$ is analytic in $x$ at least for $x \in B$.  Moreover, the only dependence on $v$ occurs through the definition of the function $A_v$, which is analytic in $x$ uniformly in $v$, since $\sigma(x,v)$ is constant in $x$ and continuous in $v$.  This shows that $\langle (T^{-1}K)^{j} Jf \rangle$ is analytic in $x$ for each $j$, and since \eqref{FiniteUExpansion} gives $\langle u \rangle$ as the finite sum of such terms, this proves that $\langle u \rangle$ is analytic as well, concluding the proof.

\end{proof}

Theorem \ref{MainTheorem} now follows from Proposition \ref{AttenuatedPartialXRay} and Lemma \ref{AnalyticRTE}.

\begin{proof}[Proof of Theorem \ref{MainTheorem}]
As in the proof of Proposition \ref{ConvexProp}, we begin by restricting our attention to the intersection of a single plane with $\Omega$.  Let $P$ be a plane intersecting $E$ such that $\Omega \cap P$ is convex.  

Now we can follow the proof of Proposition \ref{ConvexProp}. Note that at $E$, $Jf = 0$, since $F$ is not connected to $E$ by straight lines through $\Omega$.   It follows from \eqref{BoundaryU} that at $E$,
\begin{equation}\label{NonConvexU}
u = w
\end{equation}
where $w$ solves
\begin{equation}\label{WEqn}
w = T^{-1}(Kw + KJf)
\end{equation}
on $\Omega$.  Then the proof of Proposition \ref{ConvexProp} shows us that $u|_{E_+} = 0$ implies that $Kw + KJf = 0$ on $P \cap \Omega$. Since, by hypothesis, every point in $B$ lies in such an intersection, it follows that $Kw + KJf = 0$ on $B$.  Since \eqref{BoundaryU} says that $u = w + Jf$, it follows that $Ku = 0$ on $B$.

Now suppose $a \in \partial A \cap B$, and pick $\theta$ and $v$ so that $(a,\theta,v) \in \partial A_-$.  Using \eqref{BoundaryU} once more, we get
\begin{align*}
u(a,\theta,v) &= w(a,\theta,v) + Jf(a,\theta,v) \\
                 &= T^{-1}(Kw + KJf)(a,\theta,v) + Jf(a,\theta,v) \\
		 &= \int_0^{\|a - a_{\theta-}\|} \alpha_v(a, a - t\theta)(Kw + KJf)(a - t\theta,v)dt + Jf(a,\theta,v). \\
\end{align*}
Since $A$ is convex, and $(a,\theta,v) \in \partial A_-$, it follows that $a - t\theta \notin A$ for all $t$. Since $Kw + KJf = 0$ on $B$, this implies that 
\[
u(a,\theta,v) = Jf(a,\theta,v) = \alpha_v(a,a_{\theta-})f(a_{\theta-},v)
\]
for $a$ in $\partial A \cap B$ and $(a,\theta, v) \in \partial A_-$. But by the same reasoning, as above, $a_{\theta-} \notin A$ and so $u(a,\theta,v) = 0$ for all $(a,\theta,v) \in \partial A_-$ with $a \in B$.   

Therefore if we restrict $u$ to $A$, we find that the hypotheses of Lemma \ref{AnalyticRTE} apply, and hence $Ku$ is analytic on $A$.  But $Ku= 0$ on $A \cap B$, which has nonempty interior, and therefore $Ku = 0$ on all of $A$. Since we already know that $Ku= 0$ on $B$, we have $Ku = 0$ on all of $\Omega$.

Using \eqref{BoundaryU}, we get that $Kw + KJf = 0$ on $\Omega$, and so \eqref{WEqn} tells us that $w = 0$ on $\Omega$.  This implies that $KJf = 0$ on $\Omega$.  Since $k >0$ on $\Omega$, this in turn implies that 
\[
\int_{S^{n-1}} \int_0^{\infty} Jf(x, \theta, v')q(v - v') \, dv'  \, d\theta = 0
\]
on $\Omega$.  By the last part of the Raman assumption, we can obtain 
\[
\langle Jf \rangle(x,v_0) = \int_{S^{n-1}} Jf(x, \theta, v_0) d\theta = 0
\]
for all $x \in \Omega$.

The last part of the proof is to show that this implies that $f(x,v_0) = 0$ for $x \in F$.  To do this, note that  \eqref{JConvo} shows that $\langle Jf \rangle(x,v_0)$ is the sum of integrals of kernels $A_i$ defined by \eqref{ADefn} against the distributions $\Phi_i$, and so it can be extended to all $x \in \mathbb{R}^n \setminus F$. Then Lemma \ref{AnalyticIntegral} shows that $\langle Jf \rangle(x,v_0)$ is analytic for $x \in \mathbb{R}^n \setminus F$, and so it follows from the fact that $\langle Jf \rangle(x,v_0) = 0$ on $\Omega$ and the fact that $x \in \mathbb{R}^n \setminus F$ is connected, that $\langle Jf \rangle(x,v_0) = 0$ for all  $x\in \mathbb{R}^n \setminus F$.  

Now \eqref{JConvo1} tells us that 
\[
\langle Jf\rangle(x,v_0)  = \int_{\partial \Omega} \alpha_{v_0}(x,y)f(y,v_0)\frac{\vec{n}(y)\cdot(y-x)}{|y-x|^n} \, dS(y).
\]
Since $\sigma(x,v_0)$ is constant in $x$, this means that 
\[
\alpha_{v_0}(x,y) = e^{-\sigma_{v_0}|x-y|} = 1 + b_{v_0}(|x-y|)|x-y|,
\]
for some differentiable function $b_{v_0}:\mathbb{R} \rightarrow \mathbb{R}$. Therefore the expression for $\langle Jf \rangle(x,v_0)$ breaks up as 
\begin{align*}
\langle Jf\rangle(x,v_0)  &= \int_{\partial \Omega} f(y,v_0)\frac{\vec{n}(y)\cdot(y-x)}{|y-x|^n} \, dS(y) \\
& \, \quad \, + \int_{\partial \Omega} b_{v_0}(|x-y|)f(y,v_0)\frac{\vec{n}(y)\cdot(y-x)}{|y-x|^{n-1}} \, dS(y).
\end{align*}
Note that these two integrals take the form of a double-layer and a (modified) single-layer potential in $\mathbb{R}^n$, respectively \cite{ColKre}.  If we set $y_0 \in \partial \Omega$ to be on a smooth part of the boundary and take the limit as $x \rightarrow y_0$ from both the interior and exterior sides of the boundary, and subtract the two results, standard results for these potentials (see e.g. theorems 2.13 and 2.7 in \cite{ColKre}) imply that we therefore obtain 
\[
\lim_{x\rightarrow y_0^+} \langle Jf \rangle(x,v_0) - \lim_{x\rightarrow y_0^-} \langle Jf \rangle(x,v_0)  =  f(y_0,v_0).
\]
But both limits on the right hand side are zero, since $\langle Jf \rangle(x,v_0)$ vanishes for $x \notin F$.  Therefore $f(y_0,v_0) = 0$ for almost all $y_0$, and the continuity of $f$ guarantees that $f(x,v_0) \equiv 0$.  This concludes the proof. 
  
\end{proof}

\section{Conditions and Conjectures}

The main result of this paper, Theorem \ref{MainTheorem}, shows that it is possible under certain conditions to be able to see around corners using observations of the scattered light. The statement of the theorem contains four major conditions, though: the corner must be ``recoverable'', the coefficients require an analyticity condition, the RTE needs a Raman-type condition to ensure termination of the collision expansion, and the domain must be restricted to a single corner.  In the final section of this paper, I want to discuss each of these conditions in slightly more detail, and speculate as to their necessity.  

\noindent \textbf{1.} Corner Type: Theorem \ref{MainTheorem} requires a piece of the ``corner'' to be recoverable from the observation set $E$.  In particular this rules out recovery around corners of the type shown in Figure 2.  The reason for this in the proof is the solution needs to be shown to be zero in an open subset of the corner using partial data from an attenuated X-ray transform, of the type described in ~\cite{Qui93c}. The necessity of this condition then probably relies on the sharpness of this result -- it is not immediately clear to this author either how to construct a counterexample or to prove a stronger result.  

\noindent \textbf{2.} Analyticity Conditions: As written $\sigma$ must be constant, which may not be optimal. On the other hand, it would seem that some sort of analyticity condition is necessary, for at least two reasons. First, ~\cite{Bom} shows that generalized X-ray transforms with local data can fail to be invertible if the weights are not analytic. Second, observations from $E$ cannot recover more than $\langle u \rangle$ in the subset of $\Omega$ directly visible from $E$; if $\langle u \rangle \equiv 0$ on that set, then $u|_{E_{+}}$ must be identically zero.  But $\langle u \rangle \equiv 0$ is equivalent to $\langle Jf \rangle \equiv 0$ (the forward implication is in the proof of Theorem \ref{MainTheorem}, and the reverse implication can be deduced from \eqref{BoundaryU}), and insofar as $\langle Jf \rangle$ is a sort of convolution of $f$, its value on some open set does not determine $f$ without some sort of an analyticity condition. This would seem to indicate that this result does not hold for non-analytic coefficients in general.  

\noindent \textbf{3.} Raman Scattering: The chief contribution of the (strict) Raman scattering assumption is to ensure that the series expression for the solution $u$ is finite, and hence the analyticity of each term shows the analyticity of the sum.  As noted above, the strict Raman scattering assumption, and its combination with isotropicity, is somewhat unphysical; it could perhaps be replaced, but the current nature of the proof seems to require some sort of condition to terminate the collision expansion in order to guarantee analyticity. 

\noindent \textbf{4.} Single Corner Condition: Can sources be recovered around multiple corners?  Note that the current proof works for a source hidden around one corner, but fails for sources hidden around multiple corners, because observations at $E$ can show that $\langle u \rangle$ is zero in the set recoverable from $E$, but $\langle u \rangle$ isn't expected to be analytic on any non-convex set containing both $F$ and any part of the recoverable set from $E$.  Nevertheless it seems odd that recovery can proceed around one corner but not a second, and not obvious how to construct a counterexample.

\end{document}